\documentclass[11pt]{article}
\usepackage{amsfonts, latexsym, amssymb,amsmath,amsthm}
\usepackage[dvips]{graphicx,color}
\usepackage{caption}
\makeatletter
\begin{document}
\title{Asymptotics of Integrals of Some Functions Related to the Degenerate Third Painlev\'e Equation}
\author{A.~V.~Kitaev and A.~Vartanian \\
Steklov Mathematical Institute, Fontanka 27, St. Petersburg 191023, Russia \\
Department of Mathematics, College of Charleston, Charleston, SC 29424, USA \!\!\!\!
\thanks{The first author is supported by the Russian Foundation for Basic Research, grant
No. 16-01-00296. E-mail: \texttt{kitaev@pdmi.ras.ru}}}
\date{November 9, 2018}
\maketitle
\begin{abstract}
\noindent
It is shown how to calculate asymptotics of integrals over the positive semi-axis of two functions
related to the Degenerate Third Painlev\'e Equation (dP3). As an example, the corresponding
results for the meromorphic solution of the dP3 vanishing at the origin are presented.
\vspace{24pt} 
\vspace{24pt} \\
{\bf Short title:} Integrals of the Degenerate Third Painlev\'e Functions
\vspace{24pt} 
\end{abstract}
\newpage
\setcounter{page}2
\section{Introduction} \label{sec:Introduction}
The degenerate third Painlev\'{e} equation can be written in the following form
\cite{KV2004,KV2010},
\begin{equation}\label{eq:dp3}
u^{\prime \prime}=\frac{(u^{\prime})^2}{u}-\frac{u^\prime}{\tau}+\frac{1}{\tau}
(-8 \epsilon u^{2}+2ab)+\frac{b^{2}}{u},
\end{equation}
where $u=u(\tau)$, the primes denote differentiation with respect to $\tau$, $a \in \mathbb C$
and $b>0$ are parameters, and $\epsilon=\pm 1$. In most instances, the $\tau$ dependencies
are suppressed; e.g., the notation $u$ connotes $u=u(\tau)$.

There is another form of Equation~\eqref{eq:dp3}, namely,
\begin{equation}\label{eq:dp3f}
b^2\tau^2\left(f^{\prime\prime}-2b^2\right)^2+\left(8f+i\epsilon b(2ai-1)\right)^2
\left((f^\prime)^2-4b^2f\right)=0,
\end{equation}
where $i^2=-1$. Equation~\eqref{eq:dp3f} coincides with the one given in \cite{KV2004} via the
re-scalings $f \to i \epsilon bf/2$ and $a \to a-i$. It occurs because of a slight difference in the
definition of the function $f$; more precisely, for any solution $u$ of Equation~\eqref{eq:dp3},
define the functions (see \cite{KV2004}, p. 1198)
\begin{equation}\label{eq:u+}
u_+ = \frac{i\epsilon b}{8u^2}\left(\tau(-u^\prime-ib)-(2ai-1)u\right)
\end{equation}
and
$$
u_-=\frac{i\epsilon b}{8u^2}\left(\tau(u^\prime-ib)-(2ai+1)u\right),
$$
which solve Equation~\eqref{eq:dp3} for $a=a_{+} := a+i$ and $a=a_{-} := a-i$, respectively.
One proves that the function $f=u_+u$ solves Equation~\eqref{eq:dp3f}, whilst the function
$i\epsilon buu_-/2$ solves the equation for the function $f$ presented on p.~1168 of
\cite{KV2004}. Conversely, suppose that $f$ is a solution of Equation~\eqref{eq:dp3f}; then,
\begin{equation}\label{eq:f-to-u}
u=\frac{f^\prime}{2ib} - \frac{\epsilon \tau (f^{\prime\prime}-2b^2)}{2(8f+i\epsilon b(2ai-1))}
\end{equation}
solves Equation~\eqref{eq:dp3}, and $u_+=u-f^{\prime}/(ib)$.

Due to the works~\cite{O2,JM}, there's another well-known class of equations that are quadratic
with respect to the second derivative, that are equivalent to the Painlev\'e equations, namely,
the so-called $\sigma$-forms of the Painlev\'e equations, which are related with their
Hamiltonian structures and the $\bf\tau$-functions. In this paper, the $\sigma$-forms
of Equation~\eqref{eq:dp3} are not discussed.

The equation that is equivalent to Equation~\eqref{eq:dp3} is specified in the work~\cite{CSc}
(see p.~75 of \cite{CSc}) as SD-III.b (5.66) under the conditions (5.68) and denoted by
CD-III.A. One learns {}from this work that Equation~\eqref{eq:dp3} was first discovered by F.
Bureau \cite{B1972} via the direct Painlev\'e analysis: he also found a relation of this equation
with Equation~\eqref{eq:dp3}, which is, without interpreting the functions $u_\pm$ as B\"acklund
transformations, equivalent to our formulae. The transformation~\eqref{eq:f-to-u} may, in fact,
be new. It should be noted that the derivation of Equation~\eqref{eq:dp3f} in \cite{KV2004} is
based on the Hamiltonian structure of Equation~\eqref{eq:dp3} and, indirectly, its isomonodromy
deformations.

The definition of the function $f$ can be re-written as
$$
f=u_+u=-\tau\frac{i\epsilon b}{8}\left(\frac{u^\prime}u-\frac1\tau+i\left(\frac{2a}
\tau+\frac{b}u\right)\right);
$$
equivalently,
\begin{equation}\label{eq:f-ln-A}
f=-\frac{i\epsilon b}{8}\tau\frac{d}{d\tau}\ln A(\tau),\qquad A(\tau) :=
\frac{u}\tau\,e^{i\varphi},\qquad \varphi^\prime := \frac{2a}\tau+\frac{b}u,
\end{equation}
where the functions $A$ and $\varphi$ are introduced in Proposition~1.2 of \cite{KV2004}
in connection with isomonodromy deformations. Integrating along a contour $\mathcal
L(\tau_0,\tau)$ connecting points $\tau_0$ to $\tau$, one arrives at, {}from the third equation
in \eqref{eq:f-ln-A}, and, after division by $\tau$, the first equation in \eqref{eq:f-ln-A},
\begin{equation}\label{eq:integrals-def}
\int_{\mathcal L(\tau_0,\tau)}\left(\frac{2a}{\tilde\tau}+\frac{b}{u(\tilde\tau)}\right)d\tilde\tau=
\varphi|_{\mathcal L(\tau_0,\tau)},\qquad
\int_{\mathcal L(\tau_0,\tau)} \frac{f(\tilde\tau)}{\tilde\tau}\,d\tilde\tau=
\frac{\epsilon b}8\left.\left(\varphi-i\ln\frac{u}{\tau}\right)\right|_{\mathcal L(\tau_0,\tau)}.
\end{equation}

The main goal of this paper is to explain how one can evaluate these integrals. Towards this
end, one has to explain how to calculate the deviation of the functions $\varphi$ and $u$
along $\mathcal L(\tau_0,\tau)$. In this paper, the aforementioned problem is considered
asymptotically, that is, when the limits of integration belong to small neighbourhoods of
the singular points, $0$ and $\infty$, of Equation~\eqref{eq:dp3}. For this purpose, one
requires asymptotics of the functions $u$ and $\varphi$.

Asymptotics of the function $u$ were studied in \cite{KV2004,KV2010}: the corresponding
asymptotics for the function $\varphi$ can also be extracted {}from these papers. In order
to do so, recall that in Proposition~1.2 of \cite{KV2004} there was one more function:
$$
B(\tau)=-\frac{u}\tau\,e^{-i\varphi};
$$
therefore, the function $\varphi$ can be presented as
\begin{equation}\label{eq:varphi-A-B}
\varphi=-\frac{i}2\ln\left(-\frac{A}{B}\right)=-i\ln \left(\frac{\sqrt{-AB}}{B} \right).
\end{equation}
The final transformation of the above equation is necessary because it is for the functions
$\sqrt{-AB}$ and $B$ that asymptotic results are given in Proposition 4.3.1, Corollary 4.3.1,
and Propositions 5.5 and 5.7 of \cite{KV2004}. It is important to note that in Appendix B of
the subsequent paper~\cite{KV2010}, inconsistencies in the paper~\cite{KV2004} were located
and rectified. Furthermore, as explained in Section~7 of \cite{K2018}, due to the discrepancy
in the definition of the canonical solutions and the corresponding linear ODE, one has to add
to the asymptotics of the function $\varphi$, obtained with the help of the results in
\cite{KV2004,KV2010}, the term $a\ln\tau$.

The integral analogous to the first one in \eqref{eq:integrals-def}, but for the second
Painlev\'e equation, was calculated in \cite{BBDI}; in the latter case, however, the analogue
of Equation~\eqref{eq:dp3f} does not exist.
\section{Meromorphic Solution Vanishing at the Origin}

In the previous section, the general scheme allowing one to calculate the
integrals~\eqref{eq:integrals-def} was presented; however, for every particular solution
and contour of integration, there are special questions that must be addressed. Here,
one simple, yet interesting, example of such a calculation is considered. Note that in
this section $\epsilon=+1$.

It is proved in \cite{K2018} that for all $a\in\mathbb{C}\setminus i\mathbb Z$, there
exists the unique odd meromorphic solution of Equation~\eqref{eq:dp3} such that
$u(0)=0$. The asymptotic calculation of the integrals for this solution is considered
by taking the simplest contour, $\mathcal L(0,\tau)=[0,\tau]$, $\tau\in\mathbb R_+$,
$\tau\to+\infty$.

Consider, first, the case $a\in\mathbb R\setminus\{0\}$. For $a=0$, the solution holomorphic
in a neighbourhood of $\tau=0$ and vanishing at $\tau=0$ does not exist. For $a>0$, such
a solution has an infinite number of poles on the real axis, which can be deduced {}from the
results of \cite{KV2010}. Therefore, only the case $a<0$ is considered below. Henceforth, by
$u(\tau)$ is meant only this special solution.

It is proved in \cite{K2018} that $u(\tau)$ has neither poles nor zeros on the real axis, except
at the origin, where, by definition, $u(0)=0$. It is easy to establish {}from Equation~\eqref{eq:dp3}
that $u(\tau)$ is real for real $\tau$, and $u^\prime(0)=-b/2a$. In this case, $b>0$ and $a<0$,
so that it is obvious that $u(\tau)>0$ for $\tau>0$ and $u(\tau)<0$ for $\tau<0$, since it is an
odd function. Using the Taylor expansion for the function $u(\tau)$ (see Equation~(23) of
\cite{K2018}), one finds that
$$
\underset{\tau\to0}\lim\left(\frac{2a}{\tau}+\frac{b}{u}\right)=0;
$$
therefore, the integral of the function $2a/\tau+b/u$ exists on the real segment $[0,\tau]$.

Since the function $u$ is real, the functions $u_\pm$ are complex conjugates, $u_+=\bar u_-$;
moreover, Equation~\eqref{eq:u+} implies that $u_+$ does not have poles on the real axis. The
function $f(\tau)/\tau$ vanishes as $\tau\to0$, since $u(0)=u_+(0)=0$. Therefore, the integral
of the function $f(\tau)/\tau$ is properly defined on the real segment $[0,\tau]$.

Now, using Equation~\eqref{eq:varphi-A-B} and Proposition~4.3.1 of \cite{KV2004} (with the
corrections indicated above), one finds that
\begin{align}\nonumber
\varphi (\tau) &\underset{\tau\to+\infty}=3b^{1/3}\tau^{2/3}+2a\ln(\tau^{2/3})-
\frac{\ln(2+\sqrt{3})}{\pi}\ln\left(1-e^{2\pi a}\right)\\
&+2a\ln{2}-a\ln(b^{1/3})+\pi+i\ln\left(g_{11}^2\left(1-e^{2\pi a}\right)\right)
+o(\tau^{-\delta}),
\label{eq:varphi-infty}
\end{align}
where $g_{11}$ is the monodromy parameter introduced in \cite{KV2004} (in this
context, it might be viewed as the constant of integration), and $\delta \! > \! 0$.
Equation~\eqref{eq:varphi-A-B} and Proposition~(5.5) of \cite{KV2004} (with the
additive correction term $a\ln\tau$) give rise to the following result:
\begin{equation}\label{eq:varphi-0}
\varphi(0)=\frac{3\pi}2-a\ln{b}+2a\ln{2}+2\mathrm{Arg}(\Gamma(1+ai))+
i\ln\left(g_{11}^2\left(1-e^{2\pi a}\right)\right),
\end{equation}
where $\Gamma(\cdot)$ is the Gamma-function \cite{BE}. Subtracting
Equation~\eqref{eq:varphi-0} {}from Equation~\eqref{eq:varphi-infty}, one arrives at
\begin{align} \label{eq:varphi-tau-varphi-0}
\int_{0}^\tau \! \left(\frac{2a}{\tilde\tau}+\frac{b}{u(\tilde\tau)}\right)d\tilde\tau&=
3b^{1/3}\tau^{2/3}+2a\ln(b^{1/3}\tau^{2/3}) - \frac{\ln(2+\sqrt{3})}{\pi} \ln
\left(1-e^{2\pi a}\right) \nonumber \\
&-\frac{\pi}2-2\mathrm{Arg}(\Gamma(1+ai))+o(\tau^{-\delta}).
\end{align}
One recalls that the function $\mathrm{Arg}(\Gamma(1+ai))$ is defined as a continuous
function of $a$, such that $\mathrm{Arg}(\Gamma(1+ai))=\arg(\Gamma(1+ai))$ for $a
\in (-3-\pi/2,0)$: when $a$ decays {}from $0$ to $-3-\pi/2$, both arguments decay {}from
$0$ to $-\pi$; for $a=-3-\pi/2$, the function $\arg$ suffers a jump discontinuity of $2\pi$
({}from $-\pi$ to $+\pi$), and then continues to decay, whilst $\mathrm{Arg}$ continues to
decay without this jump.

To calculate the second integral, one needs the following results:
\begin{align}\label{eq:ln-u-tau-0}
&\underset{\tau\to0}{\lim}\ln\left(\frac{u}{\tau}\right)=\ln{b}-\ln(-a)-\ln{2},\\
\label{eq:ln-u-tau-infty}
&\ln\left(\frac{u}{\tau}\right)\underset{\tau\to+\infty}=-\frac23\ln{\tau}+
\frac23\ln{b}-\ln{2}+O(\tau^{-1/3}).
\end{align}
With the help of these estimates and Equations~\eqref{eq:integrals-def}, one deduces that
\begin{align}\label{eq:int-f-Re}
\mathrm{Re}
\int_{0}^{\tau} \frac{f(\tilde\tau)}{\tilde\tau} \, d \tilde \tau&=\frac{b}8 \int_{0}^{\tau}
\! \left(\frac{2a}{\tilde\tau}+\frac{b}{u(\tilde\tau)}\right)d\tilde\tau,\\
\mathrm{Im}
\int_0^\tau \frac{f(\tilde\tau)}{\tilde\tau}\,d\tilde\tau&=\frac{b}8
\left(\ln(b^{1/3}\tau^{2/3})-\ln(-a)-O(\tau^{-1/3})\right). \label{eq:int-f-Im}
\end{align}
In Equation~\eqref{eq:int-f-Im} a minus sign is indicated in the $O$-estimate in order to
stress that it is exactly the same function as in Equation~\eqref{eq:ln-u-tau-infty}. The
results derived in \cite{KV2004} allow one to obtain the $O(\tau^{-1/3})$ terms in
Equations~\eqref{eq:ln-u-tau-infty} and \eqref{eq:int-f-Im} explicitly; in particular, let
$$
x=3^{1/2}b^{1/3}\tau^{2/3} \, \, (>0) \quad \text{and} \quad
q=\sqrt{\frac{-\ln(1-e^{2\pi a})}{2\pi}} \, \, (>0);
$$
then,
\begin{equation}\label{eq:uas-a<-1}
\begin{aligned}
O(\tau^{-1/3})&=-\frac{2q}{\sqrt{x}} \, \cos \left(3x+q^2\ln(3x)+\phi_0+
o(\tau^{-\delta})\right), \\
\phi_0&=a\ln(2+\sqrt{3})+q^2\ln(12)-\frac{\pi}{4}-\arg\left(\Gamma
\left(iq^2\right)\right).
\end{aligned}
\end{equation}
\section{Numerical Examples}

In this section, several features of the results obtained in the previous section are illustrated.

\begin{figure}[htbp]
\begin{center}
\includegraphics[height=2.5in,width=5in]{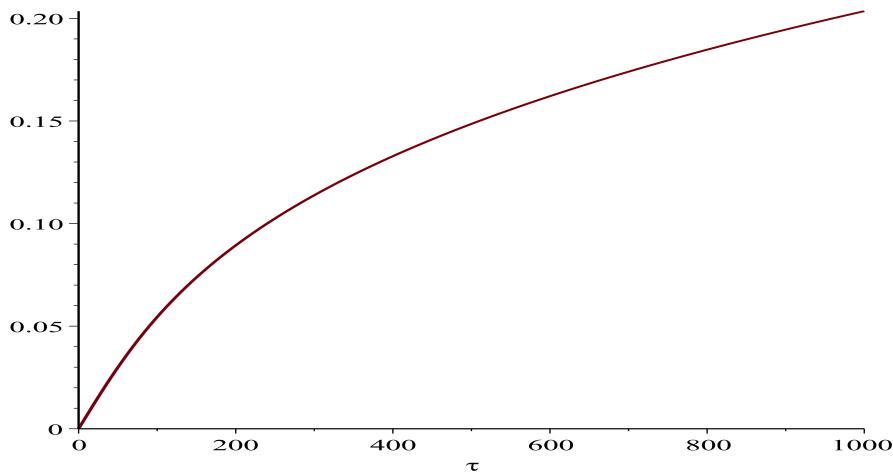}
\caption{Numerical plot of $u(\tau)$ for $a=-8$ and $b=1/100$.}
\label{fig:solam8b1d100}
\end{center}
\end{figure}

It is know \cite{KV2004} that $u(\tau)$ is, in fact, oscillating about the parabola in
Figure~\ref{fig:solam8b1d100}; however, this oscillatory structure is too fine to be
observed for $a<-1$.

\begin{figure}[htbp]
\begin{center}
\includegraphics[height=2.5in,width=5in]{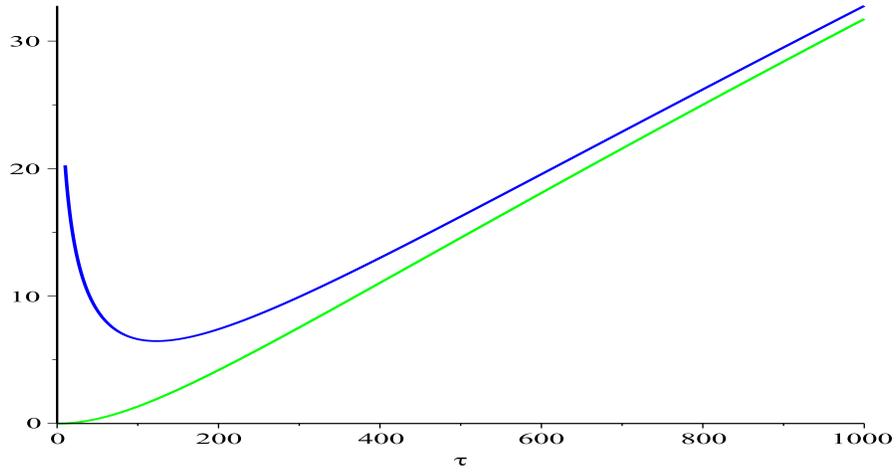}
\caption{Plot of $\int_0^\tau\left(2a/{\tilde\tau} +b/u(\tilde\tau)\right)d\tilde\tau$ for
$a=-8$ and $b=1/100$. The upper line is the asymptotics and the lower line is the
numerical plot of the integral.}
\label{fig:intam8b1d100}
\end{center}
\end{figure}

For the calculation of asymptotics, via Equation~\eqref{eq:varphi-tau-varphi-0}, in
Figure~\ref{fig:intam8b1d100} it is important to note that $\mathrm{Arg}\left(
\Gamma(1-8i)\right)=\arg\left(\Gamma(1-8i)\right)-2\pi$.

\begin{figure}[htbp]
\begin{center}
\includegraphics[height=2.5in,width=5in]{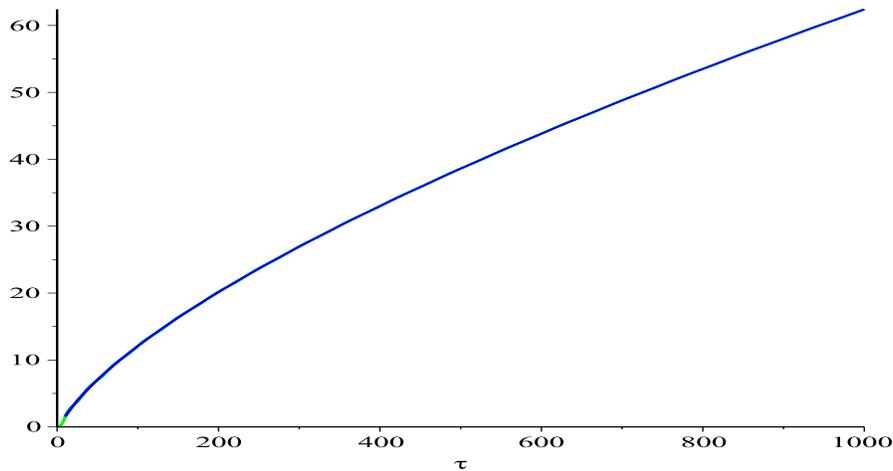}
\caption{Plot of $\int_0^\tau\left(2a/{\tilde\tau}+b/u(\tilde\tau)\right)d\tilde\tau$ for
$a=-1/8$ and $b=1/100$. The asymptotic and numerical solutions practically coincide
for $\tau>5$.}
\label{fig:intam1d8b1d100}
\end{center}
\end{figure}

According to Equation~\eqref{eq:int-f-Re}, the plots in Figures~\ref{fig:intam8b1d100}
and \ref{fig:intam1d8b1d100} for $\mathrm{Re}\int_0^\tau\left(f(\tilde\tau)/{\tilde\tau}
\right)d\tilde\tau$ coincide modulo the numeric factor $b/8$.

\begin{figure}[htbp]
\begin{center}
\includegraphics[height=2.5in,width=5in]{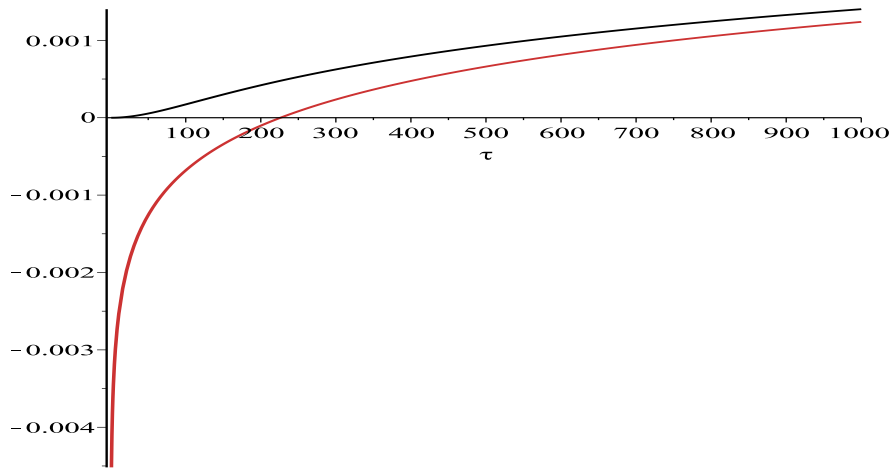}
\caption{Plot of $\mathrm{Im}\int_0^\tau\left(f(\tilde\tau)/{\tilde\tau}\right)d\tilde\tau$ for
$a=-8$ and $b=1/100$. The upper line is the numerical plot of the integral and the lower
line is the plot of its asymptotics \eqref{eq:int-f-Im}.}
\label{fig:iifam8b1d100}
\end{center}
\end{figure}

The correction term \eqref{eq:uas-a<-1}  is not observable in Figure~\ref{fig:iifam8b1d100}:
this is the general situation for all values $a<-1$; it is obviously related with the analogous
situation for $u(\tau)$.

\begin{figure}[htbp]
\begin{center}
\includegraphics[height=2.5in,width=5in]{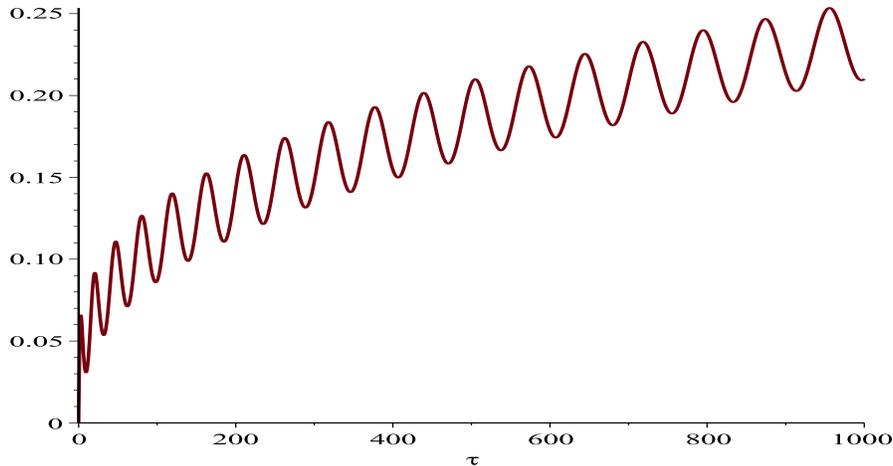}
\caption{Numerical plot of $u(\tau)$ for $a=-1/8$ and $b=1/100$.}
\label{fig:solam1d8b1d100}
\end{center}
\end{figure}

For $-1<a<0$, oscillations of the solution (which are ``hidden'' for smaller values of $a$)
are clearly seen.

Oscillations that are seen in Figure~\ref{fig:solam1d8b1d100} are not observable in
Figure~\ref{fig:intam1d8b1d100}, and, consequently, for $\mathrm{Re}\int_0^\tau
\left(f(\tilde\tau)/{\tilde\tau}\right)d\tilde\tau$.

\begin{figure}[htbp]
\begin{center}
\includegraphics[height=2.5in,width=5in]{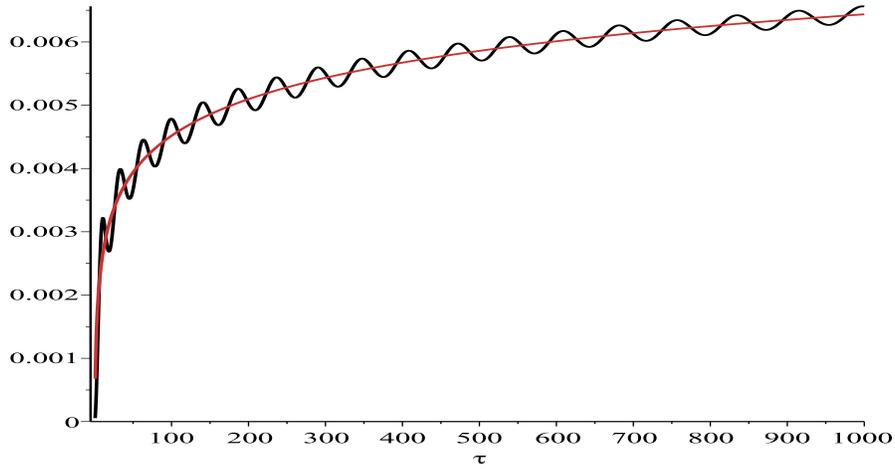}
\caption{Plot of $\mathrm{Im}\int_0^\tau\left(f(\tilde\tau)/{\tilde\tau}\right)d\tilde\tau$
for $a=-1/8$ and $b=1/100$. One sees oscillation of the numerical solution about the
asymptotic line \eqref{eq:int-f-Im} without the correction term \eqref{eq:uas-a<-1}.}
\label{fig:iifam1d8b1d100}
\end{center}
\end{figure}
\begin{figure}[htbp]
\begin{center}
\includegraphics[height=2.5in,width=5in]{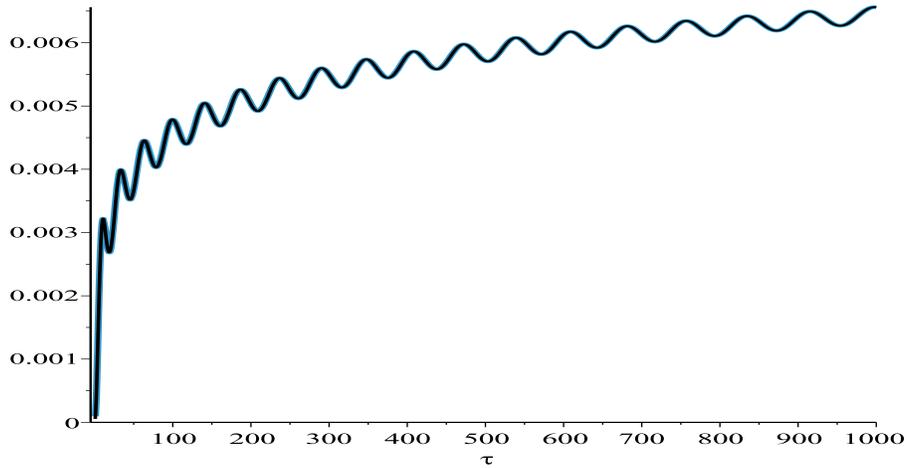}
\caption{Plot of $\mathrm{Im}\int_0^\tau\left(f(\tilde\tau)/{\tilde\tau}\right)d\tilde\tau$
for $a=-1/8$ and $b=1/100$. One notes that the numerical solution practically coincides
with the asymptotic solution with the correction term \eqref{eq:uas-a<-1}.}
\label{fig:iifcam1d8b1d100}
\end{center}
\end{figure}

Figures~\ref{fig:iifam1d8b1d100} and \ref{fig:iifcam1d8b1d100} illustrate, for $-1<a<0$,
the importance of the correction term \eqref{eq:uas-a<-1} for  $\mathrm{Im}\int_0^\tau
\left(f(\tilde\tau)/{\tilde\tau}\right)d\tilde\tau$.

Note that the value of the parameter $b>0$ is not important for observing the oscillations;
for larger values of $b$, the oscillations become faster. The value $b=1/100$ is chosen
only for the purpose of obtaining clearer figures.


\clearpage
\begin{thebibliography}{100}
\bibitem{KV2004} A.~V.~Kitaev and A.~H.~Vartanian, Connection formulae for
asymptotics of solutions of the degenerate third Painlev\'{e} equation: I,
{\it Inverse Problems}~{\bf20}, no.~4, 1165--1206 (2004).
\bibitem{KV2010} A. V. Kitaev and A. Vartanian, Connection formulae for asymptotics
of solutions of the degenerate third Painlev\'e equation: II, {\it Inverse Problems}~{\bf26},
no.~10, 105010 (2010).
\bibitem{O2} K.~Okamoto, On the $\tau$-function of the Painlev\'e equations,
{\it Physica} {\bf 2D}, 525-535 (1981).
\bibitem{JM} M.~Jimbo and T.~Miwa, Monodromy preserving deformation of linear ordinary
differential equations with rational coefficients II, {\it Physica} {\bf2D}, 407--448 (1981).
\bibitem{CSc} C. M. Cosgrove and G. Scofus, Painlev\'{e} Classification of a Class of Differential
Equations of the Second Order and Second Degree, {\it Stud. Appl. Math.} {\bf88}, no.~1,
25–87 (1993).
\bibitem{B1972} F. Bureau, \'Equations diff\'erentielles du second ordre en Y et du second
degr\'e en \"Y dont l'int\'egrale g\'en\'erale est \`a points critiques fixes, {\it Ann. Mat.
Pura Appl. (4)} {\bf91}, 163-281 (1972).
\bibitem{K2018} A.~V.~Kitaev, Meromorphic Solution of the Degenerate Third Painlev\'{e}
Equation Vanishing at the Origin, \texttt{https://arxiv.org/abs/1809.00122}.
\bibitem{BBDI} J. Baik, R. Buckingham, J. DiFranco, and A. Its, Total integrals of global
solutions to Painlev\'{e} II, {\it Nonlinearity} {\bf22}, no. 5, 1021--1061 (2009).
\bibitem{BE} A. Erdélyi, W. Magnus, F. Oberhettinger, F. G. Tricomi, \emph{Higher Transcendental
Functions}, Vol.~1, McGraw-Hill Book Company, Inc., New York-Toronto-London (1953)
(based, in part, on notes left by Harry Bateman).
\end{thebibliography}
\end{document}